\documentclass[12pt]{amsart}
\usepackage{amsmath}
\usepackage{amssymb} 
\usepackage{amscd}
\usepackage[mathscr]{eucal} 
\usepackage{latexsym}

\newtheorem{thm}{Theorem}[section]

\theoremstyle{definition}

\newtheorem{rem}[thm]{Remark}

\theoremstyle{remark}

\numberwithin{equation}{section}

\def\Hom{\text{\rm Hom}}

\def\lim{\text{\rm lim}}

\begin{document} 

\title[Algebraic cycles and the triangulated category of mixed motives]
{Algebraic cycles and the triangulated category of mixed motives}

\author[Kazuma Morita]{Kazuma Morita}
\address{Department of Mathematics, Hokkaido University, Sapporo 060-0810, Japan}
\email{morita@math.sci.hokudai.ac.jp}

\date{\today}
\maketitle
\begin{abstract}In this paper, we shall give a candidate for the $t$-structure on the triangulated category of mixed motives due to Voevodsky. 
\end{abstract}
\section{Introduction}
Let $k$ be a field of characteristic $0$ and $\mathscr{L}(k)$ denote the category of the smooth schemes over $k$. Voevodsky constructs the triangulated category of mixed motives with coefficients in $\mathbb{Q}$ (denoted by  $DM_{\text{gm}}(k)_{\mathbb{Q}}$) such that we have a canonical functor 
$$M:\mathscr{L}(k)\rightarrow DM_{\text{gm}}(k)_{\mathbb{Q}}, \quad X\mapsto M(X).$$
It has been believed that a good $t$-structure on $DM_{\text{gm}}(k)_{\mathbb{Q}}$ would capture the mixed motives of Grothendieck and this is one of the most significant problems in the field of algebraic cycles. 
\begin{rem}Although other triangulated categories of mixed motives are constructed ([H], [L]), we shall use the category of Voevodsky since it is widely spread and more accessible than others.
\end{rem}
\section{Review of $t$-structures}
For a triangulated category $\mathbb{D}$, let $\mathbb{D}^{\leq 0}$ and $\mathbb{D}^{\geq 0}$ be full subcategories of $\mathbb{D}$. We say that the couple $(\mathbb{D}^{\leq 0}, \mathbb{D}^{\geq 0})$ is a $t$-structure on $\mathbb{D}$ if the following conditions are satisfied:
\begin{enumerate}
\item $\mathbb{D}^{\leq 0}[1]\subset \mathbb{D}^{\leq 0}$ and $\mathbb{D}^{\geq 0}[-1]\subset \mathbb{D}^{\geq 0}$
\item $\Hom_{\mathbb{D}}(X,Y)=0$ for $X\in \mathbb{D}^{\leq 0}$ and $Y\in \mathbb{D}^{\geq 0}[-1]$
\item For any $X\in \mathbb{D}$, there exists a distinguished triangle $X_{0}\rightarrow X\rightarrow X_{1}\rightarrow$ in $\mathbb{D}$ such that we have $X_{0}\in \mathbb{D}^{\leq 0}$ and $X_{1}\in \mathbb{D}^{\geq 0}[-1]$.
\end{enumerate}
The full subcategory $\mathscr{A}=\mathbb{D}^{\leq 0}\cap \mathbb{D}^{\geq 0}$ is called the heart of the $t$-structure. In the next section, we shall give a  $t$-structure on $\mathbb{D}=DM_{\text{gm}}(k)_{\mathbb{Q}}$ which reflects the intersections of algebraic cycles and we expect that the heart of this $t$-structure is the mixed motives which Grothendieck has dreamed of. 
\section{Definition of a $t$-structure}
Keep the notation as in the Introduction and denote $\mathbb{D}=DM_{\text{gm}}(k)_{\mathbb{Q}}$ for simplicity. 
Let us define a $t$-structure $(\mathbb{D}^{\leq 0}, \mathbb{D}^{\geq 0})$ on $\mathbb{D}$ as follows.

$\underline{\text{Preliminary}}$

Let $\phi:A\rightarrow M(X)$ and $\psi:M(Y)\rightarrow B$ be two morphisms in $\mathbb{D}$. If we convert the orientations of $X$ and $Y$, we denote the corresponding morphisms by  $\phi^{X}:A\rightarrow M(X)$ and $\psi_{Y}:M(Y)\rightarrow B$.

$\underline{\text{Definition of }\mathbb{D}^{\leq 0}}$

The full subcategory $\mathbb{D}^{\leq 0}$ is consisted of objects $A$ of $\mathbb{D}$ such that we have
$$\phi+\phi^{X}=0:A\rightarrow M(X)[i]\qquad (\exists X\in \mathscr{L}(k), \ \exists i\in\mathbb{Z}_{\geq 0}).$$

$\underline{\text{Definition of }\mathbb{D}^{\geq 0}}$

The full subcategory $\mathbb{D}^{\geq 0}$ is consisted of objects $B$ of $\mathbb{D}$ such that  we have
$$\psi-\psi_{Y}=0:M(Y)[j]\rightarrow B\qquad (\forall Y\in\mathscr{L}(k),\ \forall j\in\mathbb{Z}_{<0}).$$
 
{\bf Acknowledgments.}
The author would like to thank Professor Masanori Asakura, Kazuya Kato and Iku Nakamura for the steadfast kindness. This research is partially supported by Grant-in-Aid for Young Scientists (B).

\end{document}